\documentclass[11pt,english]{smfart}
\usepackage{amsmath, amsthm, amsopn, amsfonts, amssymb}
\usepackage{graphicx}
\usepackage{graphpap}

\usepackage{ifpdf}
\ifpdf
\usepackage{graphicx}  \DeclareGraphicsExtensions{.pdf,.mps}

\else
\usepackage{amsfonts}
\usepackage{graphicx}  \DeclareGraphicsExtensions{.ps,.eps}

\fi

\def\pasdegrille{\let\grille = \pasgrille}
\def\ecriture#1#2{\setbox1=\hbox{#1}
\dimen1= \wd1
\dimen2=\ht1
\dimen3=\dp1
\grille #2 \box1 }
\def\aat#1#2#3{
\divide \dimen1 by 48
\dimen3=\dimen1
\multiply \dimen1 by #1
\advance \dimen1 by -\dimen3
\divide \dimen1 by 101
\multiply \dimen1 by 100
\divide \dimen2 by \count11
\multiply \dimen2 by #2 
\setbox0=\hbox{#3}\ht0=0pt\dp0=0pt
  \rlap{\kern\dimen1 \vbox to0pt{\kern-\dimen2\box0\vss}}\dimen1= \wd1
\dimen2=\ht1}
\def\pasgrille{
\count12= \dimen1 
\divide \count12 by 50
\divide \dimen2 by \count12
\count11 =\dimen2
\ 
\divide \dimen1 by 48
\setlength{\unitlength}{\dimen1}
\smash{\rlap{\ }}
\dimen1= \wd1
\dimen2=\ht1
}
\def\grille{
\count12= \dimen1 
\divide \count12 by 50
\divide \dimen2 by \count12
\count11 =\dimen2
\ 
\divide \dimen1 by 48
\setlength{\unitlength}{\dimen1}
\smash{\rlap{\graphpaper[1](0,0)(50, \count11)}}
\dimen1= \wd1
\dimen2=\ht1
}

\pasdegrille
\usepackage{epic, eepic}

\theoremstyle{plain}

\newtheorem{thm}{Theorem}
\newtheorem{prop}{Proposition}[section]
\newtheorem{cor}[prop]{Corollary}
\newtheorem{lem}[prop]{Lemma}

\newtheorem*{thmcite}{Theorem A}
\theoremstyle{definition}

\newtheorem{rem}{Remark}

\numberwithin{equation}{section}

\title[Global existence]{Global existence for energy
critical waves in $3$-d domains}
\author[N. Burq]{Nicolas Burq}
\address{Universit{\'e} Paris Sud,
Math{\'e}matiques,
B{\^a}t 425, 91405 Orsay Cedex, France et Institut Universitaire de France}
\email{Nicolas.burq@math.u-psud.fr}
\author[G. Lebeau]{Gilles Lebeau}
\address{Laboratoire J.-A. Dieudonn\'e
UMR 6621 du C.N.R.S, Universit\'e de Nice - Sophia Antipolis, Parc Valrose
06108 Nice Cedex 02, France and Institut Universitaire de France}
\email{lebeau@math.unice.fr}
\author[F. Planchon]{Fabrice Planchon}
\address{ Laboratoire Analyse, G\'eom\'etrie
  \& Applications, UMR 7539, Institut Galil\'ee, Universit\'e Paris
  13, 99 avenue J.B. Cl\'ement, F-93430 Villetaneuse}
\email{fab@math.univ-paris13.fr}

\usepackage{amssymb} \usepackage{amsmath, amsthm, amsopn, amsfonts}

\begin{document}    
\begin{abstract} We prove that the defocusing quintic wave equation,
  with Dirichlet boundary conditions, is globally well posed on
  $H^1_0(\Omega) \times L^2( \Omega)$ for any smooth (compact) domain
  $\Omega \subset \mathbb{R}^3$. The main ingredient in the proof is
  an $L^5$ spectral projector estimate, obtained recently by Smith and Sogge~\cite{SmSo06}, combined with a precise study of the boundary value problem.
\end{abstract}
   
\maketitle
   
\section{Introduction}   
\label{in}
Let $\Omega \in \mathbb{R}^3$ be a smooth bounded domain with
boundary $\partial \Omega$ and $\Delta_D$ the Laplacian acting on
functions with Dirichlet boundary conditions. We are interested in
describing the relationship between certain $L^p$ estimates for the associated spectral projector,
obtained recently by H. Smith and C. Sogge~\cite{SmSo06}, and
Strichartz inequalities for solutions to the wave equation in $\Omega$. This relationship
turns out to be very simple, natural and optimal (at least in some
range of indexes), and is closely related to an earlier remark of Mockenhaupt, Seeger and Sogge
regarding Fourier integral operators (\cite{MSS93}, Corollary 3.3). As an application we consider the critical semi-linear wave equation
(with real initial data) in $\Omega $,
\begin{equation}\label{eq.NLW}
\begin{aligned}
(\partial_t^2- \Delta) u + u^{5}&=0, \qquad \text{ in } \mathbb{R}_t \times \Omega\\
u \mid_{t=0} = u_0, \qquad \partial_t u\mid_{t=0} &=u_1 , \qquad u
\mid_{\mathbb{R}_t\times \partial \Omega} =0 ,
\end{aligned}
\end{equation}
which enjoys the conservation of energy
$$ 
E(u)(t) = \int_{\Omega}\Bigl(\frac{ |\nabla u|^2 (t,x) + |\partial_t u|^2(t,x)} 2 + \frac {|u|^{6}(t,x)} {6}\Bigr)dx = E(u) (0).
$$
Our main result reads:
\begin{thm}\label{th.1}
For any $(u_0,u_1) \in H^1_0( \Omega)\times L^2(\Omega)$ there exists
a unique (global in time) solution  $u$ to~\eqref{eq.NLW} in the space
$$ X= C^0( \mathbb{R}_t; H^1_0(\Omega))\cap C^1(\mathbb{R}_t; L^2( \Omega)) \cap L^5_{\text{loc}}(\mathbb{R};L^{10}( \Omega)). $$
\end{thm}
\begin{rem} To our knowledge, the fact that weakened dispersion
  estimates can still imply optimal (and scale invariant) Strichartz
  estimates for the solution to the wave equation was first noticed by
  the second author~\cite{Leproc06}. Observe that the results obtained
  in \cite{Leproc06}, though restricted only to the interior of strictly convex
  domains, are far more precise than the results presented here, and
  apply to the critical non linear wave equation in higher dimension
  as well.
\end{rem}
\begin{rem} The difficulty in proving Theorem~\ref{th.1} is that we
  cannot afford any loss in the Strichartz estimates we prove:
  Strichartz estimates obtained by Tataru~\cite{Ta04} for Lipschitz
  metrics (see also the results by Anton~\cite{An05} in the
  Schr\"odinger context) would certainly
  improve the classical result (well posedness for the cubic non
  linear wave equation) but it would not be enough to deal with the
  quintic nonlinearity.
\end{rem}
\begin{rem} In other dimensions, we can still apply Smith and Sogge's
  spectral projectors results, leading to other Strichartz type
  estimates. The strategy also works for obtaining semi-classical
  Strichartz-type estimates for the Schr\"odinger equation (see Burq, G\'erard and Tzvetkov~\cite{BuGeTz01-1}). These questions will be addressed elsewhere.
\end{rem}
\begin{rem} Using the material in this paper, it is rather standard to prove existence of global smooth solutions, for smooth initial data satisfying compatibility conditions (see ~\cite{SmSo95}). Furthermore, the arguments developed in this paper apply equally well to more general defocusing non linearities $f(u)= V'(u)$ satisfying
$$
|f(u) | \leq C (1+ |u|^5), \qquad |f'(u)|\leq C (1+ |u|)^4. 
$$
Finally, let us remark that our results can be localized
  (in space) and consequently hold also in the exterior of any obstacle, and we extend in this framework previous results obtained by Smith and Sogge~\cite{SmSo95} for convex obstacles.
\end{rem}
We shall denote in the remaining of this paper, for $s\geq 0$, by $H^s_D( \Omega)$ the domain of $(-\Delta_D)^{s/2}$ ($H^s_D= H^s_0( \Omega)$ for $0\leq s < 3/2$).
\par \noindent
{\bf Acknowledgments:} We thank P. G\'erard for various enlightenments about the critical wave equation.
\section{Local existence}
The local (in time) existence result for~\eqref{eq.NLW} is in fact an
easy consequence of some recent work by Smith and Sogge~\cite{SmSo06}
on the spectral projector defined by $\Pi_\lambda= 1_{\sqrt{-\Delta_D}\in [\lambda, \lambda+1[}$.
\begin{thmcite}[Smith-Sogge~\protect{\cite[Theorem 7.1]{SmSo06}}]\label{th.SmSo} 
Let $\Omega \in \mathbb{R}^3$ be a smooth bounded domain, then
\begin{equation}\label{eq.proj}
 \|\Pi_\lambda u\|_{L^5(\Omega)}\leq \lambda^{\frac 2 5} \|u\|_{L^2(
 \Omega)}.
\end{equation}
\end{thmcite}
We now derive from this result some Strichartz estimates, which are
optimal w.r.t scaling.
\begin{thm}\label{th.2}
Assume that for some $2 \leq q < + \infty$, the spectral projector $\Pi_\lambda$ satisfies
\begin{equation}\label{eq.estproj}
 \|\Pi_\lambda u\|_{L^q(\Omega)}\leq \lambda^{\delta} \|u\|_{L^(
 \Omega)}.
\end{equation}
Then the solution to the wave equation $v(t,x) = e^{it \sqrt{-
\Delta_D}}u_0$ satisfies
$$
 \|v\|_{L^q((0,2\pi)_t \times \Omega_x)}\leq C \|u_0 \|_{H_D^{\delta+\frac 1 2 - \frac 1 q}}.
$$
\end{thm}
\begin{proof} 
The proof is rather simple. In fact, one can recast the spectral
projector estimate as a square function (in time) estimate for the
wave equation, and use Sobolev in time, as was already observed in
\cite{MSS93} in  the context of variable coefficients wave equations
(note that such square function estimates are useful in their own
right in a nonlinear context: see e.g. \cite{TaAp96}, Appendix B, \cite{TaoYM}).

 Let $(e_\lambda(x))$ be the eigenbasis of $L^2( \Omega)$ consisting in eigenfunctions of $-\Delta_D$ associated to the eigenvalues $(\lambda^2)$. Let us define on this basis an abstract self adjoint operator 
$$
A (e_\lambda) = [\lambda] e_\lambda
$$
(where $[\lambda]$ is the integer part of $\lambda$). Now we prove the estimate with $v$ replaced by $\widetilde {v} = e^{it A} u_0$.  We decompose $ \widetilde v(t,x) = \sum_{k\in \mathbb{N}} v_k(t,x)$ with 
$$
\widetilde v_k(t,x)= \sum _{\lambda \in \sigma(\sqrt{-\Delta_D})\cap [k, k+1)} e^{it  k}u_{\lambda} e_\lambda(x), \qquad u_0 = \sum_\lambda u_\lambda e_\lambda (x).
$$ 
Using Plancherel formula (for $x$ fixed),
$$
\|\widetilde v(\cdot,x)\|^2_{H^s(0,2\pi)} =\sum_{k\in
  \mathbb{N}}(1+ k)^{2s}\|\widetilde{v}_k(\cdot,x)\|^2_{L^2(0,2\pi)}\,; 
$$
and consequently, using Sobolev injection in the time variable for the
first inequality, \eqref{eq.estproj} in the last inequality, with $s=\frac 1 2-\frac 1 q$, and the fact that $q\geq 2$ from line 3 to line 4,
\begin{align*}
\|\widetilde v\|^2_{L^q(\Omega; L^q(0,2\pi))} & \leq C\|\widetilde v\|^2_{L^q(\Omega; H^{s}(0,2\pi))}=\|\|\widetilde v(\cdot,x)\|^2_{H^s(0,2\pi)}\|_{L^{q/2}(M)}\\
 & \leq C \|\sum_k (1+k)^{1- \frac 2 q}\|\widetilde v_k(\cdot,x)\|^2_{L^2(0,2\pi)}\|_{L^{q/2}(\Omega)}\\
 & \leq C \sum_k (1+|k|)^{1- \frac 2 q}\|\widetilde v_k\|^2_{L^{q} (\Omega;L_t^2(0,2\pi))}\\
 & \leq C \sum_k (1+|k|)^{1- \frac 2 q}\|\widetilde v_k\|^2_{L^2((0,2\pi); L^q(\Omega))}\\
 &\leq C\sum_k \sum_{\lambda \in \sigma(\sqrt{-\Delta_D})\smash{\cap [k, k+1)}} (1+|k|)^{2 \delta+1- \frac 2 q}|u_\lambda|^2\sim \|u_0\|^2_{H^{\delta+ \frac 1 2 - \frac 1 q}_D( \Omega)}.
\end{align*}
Consequently, $e^{itA}$ is continuous from $H^{\delta+ \frac 1 2 - \frac 1 q}_D( \Omega)$ to $L^q((0,2\pi); L^q( \Omega))$. Coming back to $v= e^{it \sqrt{-\Delta_D}} u_0$, it satisfies
$$ 
(i\partial_t + A) v = (A- \sqrt{- \Delta_D}) v,\, v\mid_{t=0} =u_0,
$$ 
therefore, using Duhamel formula, Minkovski inequality and that both $e^{it \sqrt{ - \Delta_D}}$ and $(A- \sqrt{- \Delta_D})$ are bounded on  $H^{\delta+ \frac 1 2 - \frac 1 q}_D( \Omega)$,
\begin{align*}
\|v\|_{L^q(\Omega; L^q(0,2\pi))} & \leq C \|u_0\|_{H^{\delta+ \frac
    1 2 - \frac 1 q}_D( \Omega)}+ C\|(A- \sqrt{- \Delta_D})
v\|_{L^1((0,2\pi); H^{\delta+ \frac 1 2 - \frac 1 q}_D( \Omega))}\\
 & \leq
C'\|u_0\|_{H^{\delta+ \frac 1 2 - \frac 1 q}_D( \Omega)}. 
\end{align*}
\end{proof}
\begin{cor}\label{cor.1} 
Consider $u$ solution to
$$
 (\partial_t^2 - \Delta ) u =0, \quad u \mid_{\partial \Omega}=0, \quad u \mid_{t=0} = u_0, \quad \partial_t u \mid_{t=0} = u_1.
$$
Then
 \begin{equation}
 \label{eq.stri2}
 \|u\|_{L^5((0,1); L^{10}( \Omega))} \leq C \left(\|u_0\|_{H^{1}( \Omega)}+ \|u_1\|_{L^2( \Omega)}\right)\, .
\end{equation}
As a consequence, for any initial data $(u_0, u_1) \in H^1_0(\Omega) \times L^2(\Omega)$, the critical non linear wave equation~\eqref{eq.NLW} is locally well posed in 
$$
 X_T= C^0([0,T]; H^1_0( \Omega)) \cap L^5((0,T); L^{10}( \Omega)) \times
 C^0([0,T]; L^2( \Omega))$$ (globally for small norm initial data).
\end{cor}
\begin{rem} This corollary proves the uniqueness part in Theorem~\ref{th.1}
\end{rem}
To prove Corollary~\ref{cor.1}, we observe that according to Theorems A and~\ref{th.2}, the operator $\mathcal{T}= e^{\pm it \sqrt{ - \Delta_D}}$satisfies
\begin{equation}\label{eq.strichartz4bis}
\|\mathcal{T}u_0\|_{L^5((0,1) \times \Omega)} \leq C \|u_0\|_{H^{\frac 7 {10}}_D( \Omega)}.
\end{equation}
Applying the previous inequality to  $\Delta u_0$, and using the $L^p$ elliptic regularity result 
\begin{equation}
\begin{gathered}
 - \Delta u +u =f \in L^p( \Omega), \quad u \mid_{\partial \Omega} =0 \Rightarrow u \in W^{2,p}( \Omega) \cap W^{1,p}_0( \Omega)\\
\text{ and }
\|u\|_{W^{2,p}( \Omega)}\leq C \|f\|_{L^p( \Omega)},\qquad 1 <p<+\infty 
\end{gathered}
\end{equation}
we get
\begin{equation}\label{eq.strichartz4}
 \|\mathcal{T}u_0\|_{L^5((0,1);W^{2,5}( \Omega)\cap W^{1,5}_0( \Omega))} \leq C \|u_0\|_{H_D^{\frac {27} {10}}( \Omega)}
\end{equation}
and consequently by (complex) interpolation between~\eqref{eq.strichartz4bis} and~\eqref{eq.strichartz4},
\begin{equation}\label{eq.strichartz6} 
 \|\mathcal{T}u_0\|_{L^5((0,1);W^{\frac 3 {10},5}_0 (\Omega))} \leq C \|u_0\|_{H_D^{1}( \Omega)}\,;
\end{equation} 
finally, by Sobolev embedding
\begin{equation}\label{eq.strichartz8}
 \|\mathcal{T}u_0\|_{L^5((0,1);L^{10}( \Omega))} \leq C \|u_0\|_{H_D^{1}( \Omega)}.
\end{equation}
To conclude, we simply observe that
$$
 u= \cos( t\sqrt{ - \Delta_D}) u_0 + \frac { \sin(t\sqrt{ - \Delta})}{ \sqrt{ - \Delta}} u_1
$$
and $1/\sqrt{- \Delta_D}$ is an isometry from $L^2( \Omega)$ to $H^1_0( \Omega)$, leading to~\eqref{eq.stri2}. The remaining of Corollary~\ref{cor.1} follows by a standard fixed point argument with $(u, \partial_t u)$ in the space $X_T$
 with a sufficiently small $T$ (depending
 on the initial data $(u_0, u_1)$). Note that this local in time
 result holds irrespective of the sign of the nonlinearity.

Finally, to obtain the global well posedness result for small initial data, it is enough to remark that if the norm of the initial data is small enough, then the fixed point can be performed in $X_{T=1}$.
Then the control of the $H^1$ norm by the energy (which is conserved along the evolution) allows to iterate this argument indefinitely leading to global existence. Note that this
 result holds also irrespective of the sign of the nonlinearity because for small $H^1$ norms, the energy always control the $H^1$ norm. 
\section{Global existence}
It turns out that our Strichartz estimates are strong enough to extend
local to {\em global} existence for arbitrary (finite energy) data, when combined with a trace estimate
and non concentration arguments.

Before going into details, let us sketch the proof. We firstly need to refine the $L^5_t; L^{10}_x$ estimate obtained above. We shall use~\eqref{eq.strichartz6} instead.

  The usual $T T^\star$ argument and Christ-Kiselev
  Lemma~\cite{ChKi01} proves the following:
  \begin{prop}
    If $u,f$ satisfy
$$
 (\partial_t^2- \Delta)u =f, \quad u\mid_{\partial \Omega} =0, \quad u \mid_{t=0}= u_0, \quad \partial_t u \mid_{t=0} = u_1
$$
then
\begin{multline}\label{eq.besov}
 \|u\|_{L^5((0,1);W^{\frac 3 {10},5}_0 ( \Omega))}+
 \|u\|_{C^0((0,1); H^1_0( \Omega))}+\|\partial_t
 u\|_{C^0((0,1); L^2( \Omega))}\\
  \leq C
 \left(\|u_0\|_{H^1_0( \Omega)} + \|u_1\|_{L^2( \Omega)}+
 \|f\|_{L^{\frac 5 4}((0,1); W^{\frac 7 {10},\frac 5 4} ( \Omega))}\right)\,.
\end{multline}
Furthermore, \eqref{eq.besov} holds (with the same constant $C$) if one replaces the time interval $(0,1)$ by any interval of length smaller than $1$.
  \end{prop}
\begin{rem} An immediate consequence of~\eqref{eq.besov} with $f=0$
  (and Minkovski inequality) is
\begin{multline}\label{eq.faible}
 \|u\|_{L^5((0,1);W^{\frac 3 {10},5}_0 ( \Omega))}+ \|u\|_{L^\infty((0,1); H^1_0( \Omega))}+\|\partial_t u\|_{L^\infty((0,1); L^2( \Omega))}\\
  \leq C
 \left(\|u_0\|_{H^1( \Omega)} + \|u_1\|_{L^2( \Omega)}+
 \|f\|_{L^{1}((0,1); L^2( \Omega))}\right),
\end{multline}
but we will need the $L^p_t$ with $p>1$ on the righthandside of
\eqref{eq.besov} later on.
\end{rem}
\begin{proof}
We have 
$$
u(t, \cdot)= \cos( t\sqrt{ - \Delta_D})u_0+ \frac{\sin( t\sqrt{ - \Delta_D})} {\sqrt{- \Delta_D}} u_1 + \int_0^t\frac{\sin( (t-s)\sqrt{ - \Delta_D})} {\sqrt{- \Delta_D}}f(s, \cdot) ds.
$$
 The contributions of $(u_0,u_1)$ are easily dealt with, as previously. Let us focus on the contribution of
$$
\int_0^t \frac{e^{i (t-s)\sqrt{ - \Delta_D}}} {\sqrt{- \Delta_D}}f(s, \cdot) ds.
$$ 
Denote by $\mathcal{T}= e^{it \sqrt{-\Delta_D}}$; interpolating between~\eqref{eq.strichartz4bis} and~\eqref{eq.strichartz4}, 
$$
\|\mathcal{T} f\|_{L^5((0,1); W^{2- \frac 7 {10},5}( \Omega)\cap W_0^{1,5}( \Omega))}\leq C \|f\|_{H^2_D( \Omega)}.
$$
 Let $u_0\in L^2$. Then there exist $v_0\in H^2_D( \Omega)$ such that 
$$
-\Delta v_0 = u_0, \qquad \|u_0\|_{L^2}\sim \|v_0\|_{H^2};
$$
as a consequence, from $\mathcal{T}u_0= \Delta \mathcal{T}v_0$,
\begin{align*}
 \|\mathcal{T}u_0\|_{L^5((0,1); W^{-\frac 7{10},5}( \Omega))} & \leq C
 \|\mathcal{T}v_0\| _{L^5((0,1);W^{2- \frac 7 {10},5}( \Omega))\cap
 W_0^{1,5}( \Omega))}\\
 & \leq C \|v_0\|_{H^{2}_D( \Omega)}\sim \|u_0\|_{L^2( \Omega)}.
\end{align*}
By duality we deduce that the operator $\mathcal{T}^*$ defined by
$$
\mathcal{T}^* f= \int_0^1 e^{-is\sqrt{- \Delta}} f(s, \cdot) ds
$$ is bounded from $L^{\frac{5}{4}}((0,1); W^{\frac 7{10},\frac{5}{4}}( \Omega))$ to
$L^2( \Omega)$ (observe that $ W^{\frac 7{10},\frac{5}{4}}( \Omega)=  W_0^{\frac
  7{10},\frac{5}{4}}( \Omega)$); using~\eqref{eq.strichartz6} and boundedness of
 $\sqrt{-\Delta_D}^{-1} $ from $L^2$ to $H^1_D$, we
obtain 
$$
\|{\mathcal{T} {(\sqrt{- \Delta_D})^{-1}}\mathcal{T}^*}
f\|_{L^{5}((0,1); W^{\frac 3{10},5}_0( \Omega))\cap
  L^\infty((0,1); H^1_0( \Omega))}\leq C\|f\|_{L^{\frac{5}{4}}((0,1); W^{\frac 7{10},\frac{5}{4}}_0(
  \Omega))},
$$
and 
$$
\|\partial_t {\mathcal{T}} (\sqrt{-\Delta_D})^{-1}\mathcal{T}^*f\|_{L^\infty((0,1); L^2( \Omega))}\leq
C\|f\|_{L^{\frac{5}{4}}((0,1); W^{\frac 7{10},\frac{5}{4}}_0(
  \Omega))}.
$$
 But
$${\mathcal{T}} (\sqrt{- \Delta_D})^{-1}\mathcal{T}^*f(s, \cdot)= \int_0^1 \frac{e^{i(t-s)\sqrt{- \Delta_D}}} {\sqrt{- \Delta_D}} f(s, \cdot) ds
$$
and an application of Christ-Kiselev lemma~\cite{ChKi01} allows to transfer this property to the operator
$$
f \mapsto \int_0^t \frac{e^{i(t-s)\sqrt{- \Delta_D}}} {\sqrt{- \Delta_D}} f(s, \cdot) ds.  $$
\end{proof}
Now we remark that if $f= u^5$, we can estimate
\begin{equation}
\begin{aligned}
\|u^5\|_{L^{\frac 5 4}((0,1);L^{\frac {30} {17}}( \Omega))}&\leq  \|u\|^4_{L^5((0,1); L^{10}( \Omega))} \|u \|_{L^\infty((0,1); L^6( \Omega))}\,,\\
\|\nabla_x (u^5)\|_{L^{\frac 5 4}((0,1);L^{\frac {10} {9}}( \Omega))}&=5\|u^4\nabla_x u\|_{L^{\frac 5 4}((0,1);L^{\frac {10} {9}}( \Omega))}\\
&\leq 5 \|u\|^4_{L^5((0,1); L^{10}( \Omega))} \|u \|_{L^\infty((0,1); H^1( \Omega))}\,.
\end{aligned}
\end{equation}
Interpolating between these two inequalities yields
\begin{equation}\label{eq.bonstrichartz} 
\|u^5\|_{L^{ \frac 5 4}((0,1); W^{\frac 7 {10},\frac 5 4}(
  \Omega))}\leq C \|u\|^4_{L^5((0,1); L^{10}( \Omega))} \|u \|^{\frac 3
  {10}}_{L^\infty((0,1); L^6( \Omega))}\|u \|^{\frac 7
  {10}}_{L^\infty((0,1); H^1( \Omega))}\,.
\end{equation}
Following ideas of Struwe~\cite{St88}, Grillakis~\cite{Gr90} and Shatah-Struwe~\cite{ShSt93, ShSt94}, we will localize these estimates on small light cones and use the fact that the $L^\infty_t;L^6_x$ norm is small in such small cones.
\begin{rem} 
In the argument above, we need~\eqref{eq.besov} whereas~\eqref{eq.faible} would not be sufficient.
\end{rem} 
\subsection{The $L^6$ estimate}
In this section we shall always consider  solutions in 
\begin{equation}\label{eq.espace}
 X_{<t_0}= C^0([0,t_0); H^1_0( \Omega)) \cap L^5_{loc}([0,t_0); L^{10}( \Omega)) \times
 C^0([0,t_0); L^2( \Omega))
\end{equation}
of~\eqref{eq.NLW} having bounded energy and
obtained as limits in this space of smooth solutions to the analog
of~\eqref{eq.NLW} where the non linearity and the initial data have been smoothed
out. Consequently all the integrations by parts we will perform will
be licit by a limiting argument. 
\subsubsection{A priori estimate for the normal derivative}
We start with an a priori estimate on finite energy solutions
of~\eqref{eq.NLW}, which is a consequence of the uniform Lopatinski condition.
\begin{prop} \label{prop.apriori}
Assume that $u$ is a weak solution to~\eqref{eq.NLW}.
Then we have
\begin{equation} \label{eq.apriori}
\Bigl\|\frac{ \partial u} {\partial \nu} \Bigr\|_{L^2((0,t_0)\times\partial \Omega)}\leq C E(u)^{1/2}
\end{equation}
where $\frac{ \partial u} {\partial\nu}$ is the trace to the boundary of the exterior normal derivative of $u$.
\end{prop} 
\begin{proof}
Take $Z\in C^\infty ( \Omega; T\Omega)$ a vector field whose restriction to $\partial\Omega$ is equal to $\frac \partial { \partial \nu}$ and compute for $0<T<t_0$
\begin{multline*}
 \int_{0}^{T} \int_{\Omega} [(\partial_t^2- \Delta), Z] u(t,x)\cdot {u} (t,x) dx dt\\
 =\int_0^{T}\int_{\Omega} \left((\partial_t^2- \Delta) Z u(t,x)- Z(\partial_t^2- \Delta) u(t,x)\right) {u} (t,x) dx dt\,.
\end{multline*}
Integrating by parts, we obtain
\begin{multline}\label{eq.2.4}
 \int_{0}^{T} \int_{\Omega} [(\partial_t^2- \Delta), Z] u(t,x) \cdot u(t,x) dx dt= \int_0^{T} \int_{\partial \Omega} \Bigl|\frac{\partial u}{\partial \nu}\Bigr|^2  d\sigma  dt\\
{}+\int_0^{T}\int_{\Omega} -(Z u)u^5(t,x)+ Z (u^5) u (t,x) dx dt+\left[\int_\Omega \partial_t (Zu)\cdot u dx \right]_0^{T}-\left[\int_\Omega  (Zu)\cdot \partial_t u dx \right]_0^{T}
.
\end{multline}
Remark now that if $Z= \sum_{j} a_{j}(x) \frac {\partial} { \partial
  x_j}$, then integration by parts yields (using the Dirichlet boundary condition)
\begin{align}\label{eq.2.1}
 \Bigl|\int_{\Omega} -(Z u)u^5(t,x)+ Z (u^5) u (t,x) dx dt\Bigr| & =  \frac 4 6 \Bigl|\int_0^{T}\int_\Omega Z(u^6)(t,x) dx \Bigr|\\
 & =  \frac 4 6\Bigl| \int_0^{T}\int_\Omega \sum_j \frac{ \partial a_j} { \partial x_j} u^6 dx \Bigr|\leq C E(u)\nonumber
\end{align}
while
\begin{equation}\label{eq.2.2}
\Bigl|\left[\int_\Omega \partial_t (Zu) u dx \right]_0^{T}-\left[\int_\Omega  (Zu) \partial_t u dx \right]_0^{T}\Bigr| \leq C E(u),
\end{equation}
and $ [(\partial_t^2- \Delta), Z]= -[ \Delta, Z]$ as a second order differential operator in the $x$ variable is continuous from $H^1_0( \Omega)$ to $H^{-1}( \Omega)$ and consequently
\begin{equation}\label{eq.2.3}
 \Bigl|\int_{t=0}^{T} \int_{\Omega} [(\partial_t^2- \Delta), Z] u(t,x) \overline{u} (t,x) dx dt\Bigr| \leq C E(u).
\end{equation}
As the constants are uniform with respect to $0<T<t_0$, collecting~\eqref{eq.2.1},~\eqref{eq.2.2}~\eqref{eq.2.3} and~\eqref{eq.2.4} yields~\eqref{eq.apriori}.
\end{proof}
\subsubsection{The flux identity} By time translation, we shall assume later that $t_0=0$.
 Let us first define
$$
 Q = \frac{|\partial_tu|^2 + |\nabla_xu|^2 } 2 + \frac {|u|^6} 6 + \partial _t u (\frac x t \cdot \nabla_x) u,
$$
\begin{equation}\label{eq.P} 
P = \frac x t \Bigl( \frac{|\partial_tu|^2 - |\nabla_xu|^2 } 2 + \frac {|u|^6} 6\Bigr)
+ \nabla_x u  \Bigl( \partial_t u + (\frac x t \cdot \nabla_x) u + \frac u t \Bigr)\,,
\end{equation}
\begin{align*}
   D_T & = \{x; |x| <-T\},\\
 K_S^T & = \{(x,t); |x|<-t, S<t<T\}\cap \Omega\\
M_S^T & =\{x;|x|=-t, S<t<T\} \\
 \partial K_S^T & =  (([S,T]\times \partial \Omega)\cap K_S^T) \cup D_T \cup D_S \cup M_S^T\\
e(u) & = \Bigl(\frac{ |\partial_t u|^2 + |\nabla_x u|^2 } 2 + \frac  {|u|^6} 6, - \partial_t u \nabla _x u\Bigr) 
\end{align*}
and the Flux across $M_S^T$
$$
\text{ Flux }(u, M_S^T)= \int_{M_S^T} \langle e(u), \nu\rangle   d\sigma(x,t)
$$
where 
$$
\nu=\frac 1{ \sqrt 2 |x|} (-t,-\frac{x} t)=\frac 1{ \sqrt 2 |x|} (|x|,\frac{x} {|x|}) 
$$ is the outward normal to $ M_S^T$ and $d\sigma(x,t)$ the induced measure on $M_S^T$. Remark that 
\begin{align}\label{eq.fluxest}
\text{ Flux }(u, M_S^T) & =  \int_{M_S^T}\frac{|\partial_t u|^2 + |\nabla_x u|^2} 2 + \frac {|u|^6} 6 - \partial_t u \frac {x} { |x|} \cdot \nabla_x u d \sigma(x,t)\\
 & = \int_{M_S^T}\frac1 2 |\frac x {|x|} \partial_t u - \nabla_x u|^2
 + \frac {|u|^6} 6 d\sigma(x,t)\geq 0 \,.\nonumber
\end{align}
An integration by parts gives (see Rauch~\cite{Ra81} or~\cite[(3.3')]{SmSo95})
\begin{multline}\label{eq.locener}
\int_{x\in \Omega, |x|<-T}\Bigl(\frac{ |\partial_t u|^2 + |\nabla_x u|^2 } 2 + \frac  {|u|^6} 6\Bigr)(x,T) dx + \text{ Flux } (u, M_S^T) \\
=\int_{x\in \Omega, |x|<-S}\Bigl(\frac{ |\partial_t u|^2 + |\nabla_x u|^2 } 2 + \frac  {|u|^6} 6\Bigr)(x,S) dx= E_{loc}(S).  
\end{multline}
This proves that $u\mid_{M_S^T}$ is bounded in $H^1(M_S^T)\cap L^6(M_S^T)$ (uniformly with respect to $T<0$) and, since $E_{loc}(S)$ beeing a non-negative non-increasing function has a limit when $S\rightarrow 0^-$, 
\begin{equation}\label{eq.fluxestbis}
\text{Flux}(u, M_S^0)=\lim_{T\rightarrow 0^-}\text{Flux}(u, M_S^T)
=\lim_{T\rightarrow 0^-}( E_{loc} (S) - E_{loc}(T))
\end{equation} exists and satisfies      
\begin{equation}\label{eq.flux}
\lim_{S\rightarrow 0^-}\text{ Flux } (u, M_S^0)=0.
\end{equation}
\subsubsection{The $L^6$ estimate}
We are now in position to prove the classical non concentration effect:
\begin{prop}
\label{noconc}
Assume that $x_0 \in \overline{\Omega}$. Then for any 
  solution $u$ to~\eqref{eq.NLW} in the space $X_{<t_0}$, we have
\begin{equation}\label{eq.nonconc}
\lim_{t\rightarrow t_0^-} \int _{x\in \Omega \cap \{ |x-x_0| < t_0 -t\}} u^6( t,x ) dx=0.
\end{equation}
\end{prop}
\begin{proof} We follow~\cite{St88, Gr90, ShSt93, ShSt94} and simply have
  to take care of the boundary terms. We can assume that $x_0\in \partial \Omega$ as otherwise these boundary terms disappear in the calculations below (which in this case are standard). Contrarily
  to~\cite{SmSo95} we cannot use any convexity assumption to obtain
  that these terms have the right sign, but we shall use Proposition~\ref{prop.apriori} to control them. Performing a space-time translation, we can assume
  $x_0=0, t_0=0$. Integrating over $K_S^T$ the identity
$$ 
0 = \text{div} _{t,x} ( tQ+ u\partial_t u, -tP) + \frac {|u|^6} 3\,,
$$ we get (see~\cite[(3.9)--
  (3.12)]{SmSo95}),
\begin{multline*} 
0 = \int_{D_T}( TQ+ u \partial_t u )(T,x) dx - \int_{D_S} (SQ+ u \partial_t u)(S,x) dx + \frac 1 {\sqrt{2}} \int_{M_S^T} ( tQ+ u \partial_t u + x \cdot P) d\sigma(x,t)\\
- \int_{((S,T) \times \partial \Omega) \cap K_S^T} \nu(x) \cdot (tP) d\sigma (x,t)+ \frac 1 3 \int_{K_S^T} u^6 dx dt\,.
\end{multline*}
Let $T\rightarrow 0^-$. Using H\"older's inequality and the conservation of energy, we get that the first term in the left tends to $0$, whereas the last term is non negative. This yields
\begin{multline} 
 - \int_{D_S} (SQ+ u \partial_t u)(S,x) dx + \frac 1 {\sqrt{2}} \int_{M_S^0} ( tQ+ u \partial_t u + x \cdot P) d\sigma(x,t)\\
\leq \int_{((S,0) \times \partial \Omega) \cap K_S^0} \nu(x) \cdot (tP) d\sigma (x,t)\,.
\end{multline}
On the other hand, by direct calculation (see~\cite[(3.11)]{SmSo95}),
\begin{multline}
\frac 1 {\sqrt{2}} \int_{M_S^0} ( tQ+ u \partial_t u + x \cdot P) d\sigma(x,t)\\
= \frac 1 {\sqrt{2}} \int_{M_S^0}\frac 1 t | t \partial_t u + x \cdot \nabla_x u +u|^2 d \sigma(x,t) + \frac 1 2 \int_{\partial D_S} u^2(S,x) d\sigma (x)
\end{multline}
and (see~\cite[(3.12)]{SmSo95})
\begin{equation} 
- \int_{D_S} (SQ+ u \partial_t u)(S,x) dx\geq - \frac 1 2 \int _{\partial D_S} u^2 d\sigma(x) - S \int_{D_S} \frac {|u|^6(S,x)} 6 dx\,.
\end{equation}
As a consequence, we obtain
\begin{multline}\label{eq.morawetz}
(-S)\int_{D_S} \frac {|u|^6(S,x)} 6 dx + \frac 1 {\sqrt{2}} \int_{M^0_S} \frac 1 t |t \partial_t u + x\cdot \nabla_x u + u |^2 d\sigma(x,t) \\
\leq \int_{((S,0)\times \partial\Omega)\cap K_S^0} \nu(x) \cdot tP d\sigma(x)dt
\end{multline}
where $\nu(x)$ is the exterior normal to $\Omega$ at point $x$ and $d\sigma$ is the surface measure on $\partial \Omega$. Taking~\eqref{eq.P} into account (and the Dirichlet boundary condition), we obtain on $\partial \Omega$ 
$$
 t\nu(x) \cdot P = \frac 1 2 (\nu(x) \cdot x) \Bigl( \frac {\partial u} { \partial \nu}\Bigr)^2  \,.
$$
However, for $x\in \partial \Omega$, given that $x_0=0\in \partial
\Omega$, we have
$$
 \frac{x} {|x|}=t+\mathcal{O}(x),\qquad \nu(x) = \nu(0)+ \mathcal{O}(x)
$$
where $t$ is a unit vector tangent to $\partial\Omega$ at $x_0=0$. Consequently, as $\nu(0) \cdot t=0$,
$$ 
\nu(x) \cdot x = \mathcal{O}( |x|^2), \text{ for $x\in \partial \Omega$}
$$  
and the right hand side in~\eqref{eq.morawetz} is bounded (using Proposition~\ref{prop.apriori}) by 
\begin{equation}
\label{eq.bord} \sup_{x\in  K_S^0} |x|^2\times \int_{(-1,0)\times \partial \Omega} \Bigl( \frac {\partial u} { \partial \nu}\Bigr)^2  d\sigma(x) dt\leq C |S|^2 E(u).
\end{equation} 
Therefore,
\begin{equation}\label{eq.morawetzbis}
\int_{D_S} \frac {|u|^6(S,x)} 6 dx  \leq  |S|E(u)+ |S|\frac 1 {\sqrt{2}|S|} \int_{M^0_S} \frac 1 {|t|} |t \partial_t u + x\cdot \nabla_x u + u |^2 d\sigma(x,t)\,;
\end{equation}
finally, by H\"older's inequality and~\eqref{eq.fluxest}, we obtain
\begin{align*}
\frac 1 {\sqrt{2}|S|} \int_{M^0_S} \frac 1 {|t|} |t \partial_t u +
x\cdot \nabla_x u + u |^2 d\sigma(x,t)  & \leq \sqrt{2} \int_{M^0_S}
\frac{|x|}{|S|} | \frac x {|x|} \partial_t u - \nabla_x u|^2
d\sigma(x,t) \\
 & {}+ \sqrt{2} \int_{M^0_S} \frac {|u|^2}{|S||t|} d
\sigma(x,t)\\
 & \leq C \text{ Flux }(u,M_S^0) 
+C\text{ Flux }(u,M_S^0)^{1/3}\,,
\end{align*}
hence,
\begin{equation}\label{eq.morawetzter}
\int_{D_S} \frac {|u|^6(S,x)} 6 dx  \leq |S|E(u)+ C\text{ Flux }(u,M_S^0) 
+C\text{ Flux }(u,M_S^0)^{1/3}
\end{equation}
for which the right hand side goes to $0$ as $S\rightarrow 0^-$ by~\eqref{eq.flux}.
Remark that in the calculations above all integrals on $K_S^0$ and $M_S^0$ have to be understood as the limits as $T\rightarrow 0^-$ of the respective integrals on $K_S^T$ and $M_S^T$ (which exist according to \eqref{eq.fluxest}, \eqref{eq.fluxestbis}).

\end{proof}
\subsection{Global existence} 
In this section we consider $u$ the unique forward maximal solution to
the Cauchy problem~\eqref{eq.NLW} in the space $X_{<t_0}$. 
Assume that $t_0 < + \infty$ and consider  a point $x_0 \in
\overline{\Omega}$; our aim is to prove that $u$ can be extended in
a neighborhood of $(x_0,t_0)$, which will imply a contradiction. We perform a space time translation and assume that $(x_0, t_0)= (0,0)$.
\subsubsection{Localizing space-time estimates}
For $t<t'\leq 0$, let us denote by 
$$
\|u\|_{(L^p; L^{q})(K_t^{t'})}= \Bigl(\int_{s=t}^{t'} \Bigl(\int_{\{|x|<-s\} \cap \Omega } |u|^{q}(s,x) dx\Bigr)^{\frac p q} ds\Bigr)^{\frac 1 p}
$$ the $L^p_tL^{q}_x$ norm on $K_t^{t'}$ (with the usual modification if $p$ or $q$ is infinite).
Our main result in this section reads
\begin{prop}
\label{prop.L5W}
For any $\varepsilon>0$, there exists $t<0$ such that 
\begin{equation}
\|u\|_{(L^5; L^{10})(K_t^{0})}<\varepsilon.
\end{equation}
\end{prop}
\begin{proof}
We start with an extension result:
\begin{lem}\label{lem.extension} For any $x_0\in \overline{\Omega}$ there exists $r_0,C>0$ such that for any $0<r<r_0$ and any $v\in H^1_0(\Omega)\cap L^p( \Omega)$, there exist 
a function $\widetilde v_r \in H^1_0( \Omega)$ (independent of the choice of $1\leq p \leq + \infty$), satisfying
\begin{equation}
\begin{gathered}
(\widetilde v_r -v) \mid_{|x-x_0|<r\cap \Omega } =0,\\
  \int_{\Omega}|\nabla \widetilde v|^2 \leq C \int_{ \Omega}|\nabla v|^2, \qquad \|\widetilde v_r\|_{L^p( \Omega)} \leq C \|v\|_{L^p( \{|x-x_0|<r\})}.
\end{gathered}
\end{equation}
In other words, we can extend functions in $H^1_0\cap L^p$ on the ball $\{|x-x_0|<r\}$ to functions in $H^1_0(\Omega)\cap L^p( \Omega)$ with uniform bounds with respect to (small) $r>0$, for the $H^1$ {\em and} the $L^p$ norms respectively.

 Furthermore, for any $u\in L^\infty((-1,0); H^1_0( \Omega))\cap L^1_{\text{loc}}((-1,0); L^p( \Omega))$, there exist 
a function $\check u \in L^\infty((-1,0); H^1_0( \Omega))\cap L^1_{\text{loc}}((-1,0); L^p( \Omega))$, satisfying
  (uniformly with respect to~$t$)
\begin{equation}\label{eq.time}
\begin{gathered}
(\check u -u) \mid_{\{|x-x_0|<-t\}\cap \Omega} =0,\\
  \int_{\Omega}|\nabla \check u|^2(t, x) + |\partial_t \check u|^2(t,x) dx \leq C \int_{\Omega}|\nabla u|^2(t,x) + |\partial_t  u|^2(t,x)\\
 \|\check u(t, \cdot)\|_{L^p( \Omega)} \leq C \|u(t, \cdot)\|_{L^p(\Omega \cap  \{|x-x_0|<-t\})}\text{ $t$-a.s. }
\end{gathered}
\end{equation}
\end{lem}
\begin{proof}
Let us first prove the first part of the lemma. Let us first assume $x_0\in \Omega$. We use the
usual reflexion extension (suitably cut off). Fix a function $\phi\in
C^\infty_0(\frac 9 {10}, \frac {11} {10})$ equal to $1$ near $1$. Let
us define $\widetilde v_r$ (in polar coordinates $(\rho, \theta)\in
(\mathbb{R}^+ \times \mathbb{S}^2)$ centered at $x_0$) by
\begin{equation}
\label{eq.prolong}
\widetilde v_r(\rho, \theta) = 
\begin{cases}
& v_r(\rho, \theta) \text{ if } \rho<r\\
& \phi\bigl(\frac \rho {r}\bigr) v( 2r-\rho, \theta) \text{ if }  \rho>r.
\end{cases}
\end{equation}
An elementary calculation shows that (uniformly with respect to $r$)
$$
 \|\widetilde{v_r}(x)\|_{L^p( \mathbb{R}^3)}\leq C \|u\|_{L^p(\{|x|<r\})}
$$
and (using Hardy inequality to control $\partial_{x}\Bigl (\phi\bigl(\frac \rho {r}\bigr) \Bigr)\times v(2r-\rho, \theta)$),
$$
\int_{\Omega}|\nabla \widetilde v|^2 \leq C \int_{\Omega}|\nabla v|^2.
$$
Since $x_0\in \Omega$, the function $\widetilde v_r$ is in $H^1_0( \Omega)$ if $r>0$ is sufficiently small. If $x_0 \in \partial \Omega$, this is no longer true and we have to take care of the boundary condition $v\mid_{\partial \Omega}=0$. To do so the idea is to take a foliation by hypersurfaces tangent to the boundary and to perform a reflexion extension {\em tangential} to this foliation.

 We consider a change of variables $\Psi$ such that (near $x_0$) 
$$
\Omega= \Psi ( \{y=(y', y_3); y_3>0\}),\qquad \Psi(0)= x_0, \qquad \Psi'(0)= \text{Id}.
$$
Writing $y=  (rz',rz_3)$, and working in polar coordinates $z'= (\rho ,\theta)\in (\mathbb{R}^+\times \mathbb{S})$, the implicit function theorem (and the assumption $\Psi'(0)= \text{Id}$) implies that there exists a smooth function $\zeta(z_3, \theta, r)$ such that  for $z_3 \in [-\frac 1 {10}, \frac 1 {10}]$,
$$ 
\|\Psi(r\rho\theta, rz_3)\| =r \Leftrightarrow \rho= \zeta(\theta,
z_3, r).
$$
We can now define an extension in a polar coordinate for $y'$, (for $\frac {y_3} r =z_3< \frac 1 {10}$) by 
\begin{equation}
\label{eq.prolongbis}
\widetilde v_r(\rho, \theta,y_3) = 
\begin{cases}
 v(\rho, \theta, y_3) &\text{ if } \rho<r\zeta(\theta, \frac{y_3} r ,r)\\
 \phi\bigl(\frac \rho {r}\bigr) v( 2r\zeta(\theta, \frac{y_3} r,r)-\rho, \theta, y_3) &\text{ if }  \rho>r\zeta(\theta, \frac{y_3} r,r).
\end{cases}
\end{equation}
Let $\phi\in C^\infty_0(-\frac 1 {11}, \frac 1 {11})$ equal to $1$ near $0$. Then the extension we consider is
$$
(1-\phi)(\frac {y_3} r)\widetilde v_{r,1}+ \phi(\frac {y_3} r)\widetilde v_{r,2}
$$
where $ \widetilde v_{r,1}$ is the extension we built in the interior case and $\widetilde v_{r,2}$ is the extension we just built (with the foliation).
\begin{figure}[ht]
\label{fig:foliation}
$$\ecriture{\includegraphics[width=5cm]{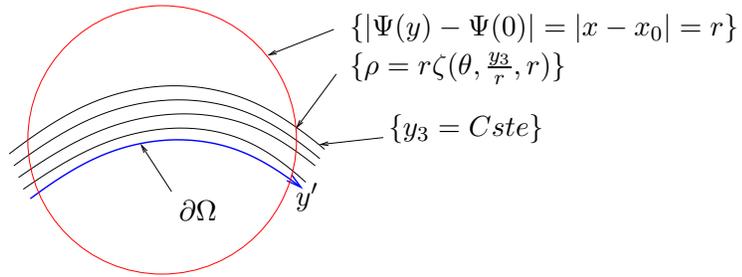}}
{\aat{45}{31}{$\{|\Psi(y)- \Psi(0)|=|x-x_0|=r\}$}\aat{45}{26}{$\{\rho=r\zeta(\theta, \frac{y_3} r,r)\}$}\aat{50}{18}{$\{y_3= Cste\}$}\aat{38}{9}{$y'$}\aat{23}{7}{$\partial \Omega$}}$$ \caption{The foliation}
\end{figure}

To obtain the second part of the lemma, we just define
$$
\check u(t, \cdot)= \widetilde u_{-t}(t,\cdot).
$$
and we can check that according to~\eqref{eq.prolong} and~\eqref{eq.prolongbis}, it satisfies~\eqref{eq.time}.
\end{proof}
Let us come back to the proof of Proposition~\ref{prop.L5W}.
 Let $\check u$ be the function given by the second part of Lemma~\ref{lem.extension}.
Then $(\check u)^5 $ is equal to $u^5$ on $K_t^{0}$ and 
\begin{equation}
\begin{aligned}
\|(\check u)^5\|_{L^{\frac 5 4}((t,t'); L^\frac{30} {17}( \Omega)}&\leq \|\check u \|^4_{ L^5((t, t');L^{10}( \Omega))}\|\check u \|_{L^\infty((t,t'); L^6(\Omega)}\\
&\leq C\|u \|^4_{ (L^5;L^{10})(K_t^{t'}))}\|u \|_{(L^\infty; L^6)(K_t^{0})}.
\end{aligned}
\end{equation}
On the other hand, $\nabla_x (\check u)^5=5 (\check u)^4 \nabla_x \check u$ and 
\begin{equation}
\begin{aligned}
\|\nabla_x(\check u)^5\|_{L^{\frac 5 4}((t,t'); L^\frac{10} {9}( \Omega)}&\leq 5\|\check u \|^4_{ L^5((t, t');L^{10}( \Omega))}\|\nabla_x\check u \|_{L^\infty((t,0); L^2(\Omega)}\\
&\leq C\| u \|^4_{ (L^5;L^{10})(K_t^{t'})}\|u\|_{L^\infty; H^1(\Omega)}.
\end{aligned}
\end{equation}
By (complex) interpolation, as in ~\eqref{eq.bonstrichartz},
$$ 
\| (\check u)^5\|_{L^{\frac 5 4}((t, {t'}); W^{\frac{7}{10}, \frac {5} 4}( \Omega))}\leq C\| u \|^4_{ (L^5;L^{10})(K_t^{t'})}\|{u}\|^{\frac 7 {10}}_{L^\infty; H^1(\Omega)}\|u \|^{\frac 3 {10}}_{(L^\infty; L^6)(K_t^{0})}.
$$
Let $w$ be the solution (which, by finite speed of propagation, coincides with $u$ on $K_t^{0}$) of 
$$
 (\partial_s^2 - \Delta )w =-(\check{u})^5,\qquad w\mid_{\partial \Omega} =0,\qquad (w-u) \mid_{s=t} =\partial_s( w-u) \mid_{s=t}=0,
$$
applying~\eqref{eq.besov}, and the Sobolev embedding $W^{\frac 3
  {10},5}(\Omega)\mapsto L^{10} ( \Omega)$, we get
\begin{multline}
\label{eq.estim}
  \|u\|_{(L^5;L^{10}) ( K_t^{t'})} \leq \|w\|_{L^5((t, t'); L^{10}(
  \Omega))} \leq C\|w\|_{ (L^5((t, t'); W^{\frac 3 {10},5})(\Omega))}\\
   \leq C E (u) +C\|u \|^4_{ (L^5; L^{10})(K_t^{t'})}\|u\|^{\frac 7 {10}}_{L^\infty; H^1(\Omega)}\|u\|^{\frac 3 {10}}_{(L^\infty; L^6)(K_t^{0})}.
\end{multline}
Finally, from Proposition \ref{noconc}, \eqref{eq.estim} and the continuity of the mapping $
 t'\in [t, 0)\rightarrow \|u\|_{ (L^5; L^{10})(K_t^{t'})} $ (which takes value $0$ for $t'=t$), there
 exists $t$ (close to $0$) such that
$$ \forall t<t'<0;  \|u\|_{ (L^5; L^{10})(K_t^{t'})}\leq 2C E (u)
$$
and passing to the limit $t' \rightarrow 0$, 
$$\|u\|_{ (L^5; L^{10})(K_t^{0})}\leq 2C E (u) .$$

As a consequence, taking $t<0$ even smaller if necessary, we obtain
\begin{equation}\label{eq.small}
  \|u\|_{ (L^5; L^{10})(K_t^{0})}\leq \varepsilon.
\end{equation}
\end{proof}
 \subsubsection{Global existence}
We are now ready to prove the global existence result.
Let $t<t_0=0$ be close to $0$ and let $v$ be the solution to the linear equation
$$
(\partial_s^2- \Delta) v=0, \qquad v\mid_{\partial \Omega} =0,\qquad (v-u)\mid_{s=t} = 0, \qquad \partial_s (v-u) \mid_{s=t}=0,
$$
then the difference $w= u-v$ satisfies 
$$
(\partial_s^2- \Delta) w=- u^5, \qquad w\mid_{\partial \Omega} =0, \qquad w\mid_{s=t} = 0, \qquad \partial_s w \mid_{s=t}=0.
$$
Let $\check u$ be the function given by Lemma~\ref{lem.extension} from ${u}$. We have 
$$
\|\check{u} \|_{L^5((t,0);L^{10} ( \Omega))}\leq C \varepsilon, \qquad \|\check u \|_{L^\infty; H^1}\leq CE(u)\,.
$$
Let $\widetilde{w}$ be the solution to
$$
(\partial_s^2- \Delta) \widetilde{w}=- \check{u}^5, \qquad \widetilde
{w}\mid_{\partial \Omega} =0, \qquad\widetilde{w}\mid_{s=t} = 0, \qquad
\partial_s \widetilde{w} \mid_{s=t}=0.
$$
By finite speed of propagation, $w$ and $\widetilde{w}$ coincide in $K_t^{0}$. On the other hand, using~\eqref{eq.faible} yields
\begin{multline}
\label{fini}
\|\widetilde{w} \|_{L^\infty((t,0);H^1)} + \|\partial_s \widetilde{w} \|_{L^\infty((t,0); L^2( \Omega))} + \|\widetilde{w} \|_{L^5((t,0); W^{\frac 3 {10},5}( \Omega))}\\
 \leq C \|\check{u}^5\|_{L^{1}((t,0); L^2( \Omega)}\leq C \|\check{u}\|^5_{L^5((t,0); L^{10} ( \Omega))}\leq C\varepsilon^5.
\end{multline}
Finally, for any ball $D$, denote by 
$$ 
E(f(s, \cdot), D)= \int_{D\cap \Omega} (|\nabla_x f|^2 + |\partial_s f|^2 + \frac {|f|^6} 3)(s,x) dx;
$$ 
since $v$ is a solution to the linear equation,
\begin{equation}
  \label{eq:lineaire}
  E(v(s, \cdot), D(x_0=0, -s)) \rightarrow 0, \qquad s \rightarrow 0^-
\end{equation}
Recalling that $u=v+\tilde w$ inside $K^0_t$, we obtain from
\eqref{eq:lineaire} and \eqref{fini} (and the Sobolev injection $H^1_0 ( \Omega \rightarrow L^6( \Omega)$) that there exists a small $s<0$
such that
$$ 
E(u(s, \cdot), D(x_0=0, -s)) < \varepsilon;
$$
but, since $(u, \partial _s u) (s, \cdot)\in H^1_0( \Omega) \times L^2( \Omega)$, we have,  by dominated convergence,
\begin{multline*}
 E(u(s, \cdot), D(x_0=0, -s))= \int_{\Omega  }1_{ \{|x-x_0|< -s\}} (x)  (|\nabla u( s, x) |^2 + |\partial_s u (s,x)|^2+ \frac {|u|^6(s, x)} 3)) dx\\
\text{ and } \lim_{ \alpha \rightarrow 0}\int_{\Omega} 1_{\{|x-x_0|<\alpha -s\} }(x) (|\nabla u( s, x) |^2 + |\partial_s u (s,x)|^2+ \frac {|u|^6(s,x)} 3)) dx\\
= \lim_{\alpha \rightarrow 0} E(u(s, \cdot), D(x_0=0, -s+\alpha))\,;
\end{multline*}
consequently, there exists $\alpha>0$ such that 
$$ 
E(u(s, \cdot), D(x_0=0, -s+\alpha))\leq 2\varepsilon.
$$
Now, according to~\eqref{eq.locener}, the $L^6$ norm of $u$ remains smaller than $2\varepsilon$ on $\{|x-x_0|<\alpha -s'\},s\leq s' <0$. As a consequence, the same proof as for Proposition~\ref{prop.L5W} shows that the $L^5; L^{10}$ norm of the solution on the truncated cone 
$$ K= \{(x,s'); |x-x_0|<\alpha-s', s<s'<0\}$$
is bounded. Since this is true for all $x_0\in \overline{\Omega}$, a compactness argument shows that 
$$ \|u\|_{L^5((s, 0); L^{10} ( \Omega))} < + \infty$$
which, by Duhamel formula shows that 
$$\lim_{s'\rightarrow 0^-} ( u, \partial_{s} u)(s', \cdot)$$ exists in $(H^1_0 ( \Omega) \times L^2( \Omega))$  and consequently $u$ can be extended for $s'>0=t_0$ small enough, using Corollary~\ref{cor.1}.
\begin{figure}[ht]
$$\ecriture{\includegraphics[width=5cm]{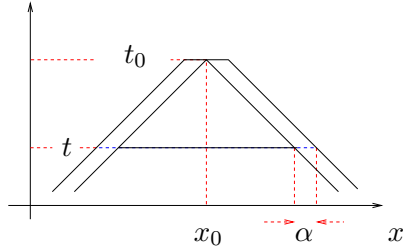}}
{\aat{25}{-2}{$x_0$}\aat{8}{9}{$t$}\aat{50}{-2}{$x$}\aat{38}{-2}{$\alpha$}\aat{16}{21}{$t_0$}}$$ \caption{The truncated cone}\label{fig.cones}
\end{figure}


\begin{thebibliography}{10}

\bibitem{An05}
Ramona Anton.
\newblock Strichartz inequalities for lipschitz metrics on manifolds and
  nonlinear {S}chr{\"o}dinger equation on domains, 2005.
\newblock preprint, {\tt arXiv:math.AP/0512639}.

\bibitem{BuGeTz01-1}
Nicolas Burq, Patrick G{\'e}rard, and Nikolay Tzvetkov.
\newblock Strichartz inequalities and the nonlinear {S}chr\"odinger equation on
  compact manifolds.
\newblock {\em Amer. J. Math.}, 126(3):569--605, 2004.

\bibitem{ChKi01}
Michael Christ and Alexander Kiselev.
\newblock Maximal functions associated to filtrations.
\newblock {\em J. Funct. Anal.}, 179(2):409--425, 2001.

\bibitem{Gr90}
Manoussos~G. Grillakis.
\newblock Regularity and asymptotic behaviour of the wave equation with a
  critical nonlinearity.
\newblock {\em Ann. of Math. (2)}, 132(3):485--509, 1990.

\bibitem{TaAp96}
Sergiu Klainerman and Matei Machedon.
\newblock Remark on {S}trichartz-type inequalities.
\newblock {\em Internat. Math. Res. Notices}, (5):201--220, 1996.
\newblock With appendices by Jean Bourgain and Daniel Tataru.

\bibitem{Leproc06}
Gilles Lebeau.
\newblock Estimation de dispersion pour les ondes dans un convexe.
\newblock In {\em Journ\'ees ``\'Equations aux D\'eriv\'ees Partielles''
  (Evian, 2006)}. 2006.

\bibitem{MSS93}
Gerd Mockenhaupt, Andreas Seeger, and Christopher~D. Sogge.
\newblock Local smoothing of {F}ourier integral operators and
  {C}arleson-{S}j\"olin estimates.
\newblock {\em J. Amer. Math. Soc.}, 6(1):65--130, 1993.

\bibitem{Ra81}
Jeffrey Rauch.
\newblock I. {T}he {$u\sp{5}$} {K}lein-{G}ordon equation. {II}. {A}nomalous
  singularities for semilinear wave equations.
\newblock In {\em Nonlinear partial differential equations and their
  applications. Coll\`ege de France Seminar, Vol. I (Paris, 1978/1979)},
  volume~53 of {\em Res. Notes in Math.}, pages 335--364. Pitman, Boston,
  Mass., 1981.

\bibitem{ShSt93}
Jalal Shatah and Michael Struwe.
\newblock Regularity results for nonlinear wave equations.
\newblock {\em Ann. of Math. (2)}, 138(3):503--518, 1993.

\bibitem{ShSt94}
Jalal Shatah and Michael Struwe.
\newblock Well-posedness in the energy space for semilinear wave equations with
  critical growth.
\newblock {\em Internat. Math. Res. Notices}, (7):303ff., approx.\ 7 pp.\
  (electronic), 1994.

\bibitem{SmSo95}
Hart~F. Smith and Christopher~D. Sogge.
\newblock On the critical semilinear wave equation outside convex obstacles.
\newblock {\em J. Amer. Math. Soc.}, 8(4):879--916, 1995.

\bibitem{SmSo06}
Hart~F. Smith and Christopher~D. Sogge.
\newblock On the $l^p$ norm of spectral clusters for compact manifolds with
  boundary, 2006.
\newblock to appear, Acta Matematica, {\tt arXiv:math.AP/0605682}.

\bibitem{St88}
Michael Struwe.
\newblock Globally regular solutions to the {$u\sp 5$} {K}lein-{G}ordon
  equation.
\newblock {\em Ann. Scuola Norm. Sup. Pisa Cl. Sci. (4)}, 15(3):495--513
  (1989), 1988.

\bibitem{TaoYM}
Terence Tao.
\newblock Local well-posedness of the {Y}ang-{M}ills equation in the temporal
  gauge below the energy norm.
\newblock {\em J. Differential Equations}, 189(2):366--382, 2003.

\bibitem{Ta04}
Daniel Tataru.
\newblock Strichartz estimates for second order hyperbolic operators with
  nonsmooth coefficients. {III}.
\newblock {\em J. Amer. Math. Soc.}, 15(2):419--442 (electronic), 2002.

\end{thebibliography}
\end{document}